\documentclass[12pt,a4paper]{amsart}

\usepackage {amsfonts}
\usepackage[english]{babel}
\def\g{\gamma}
\def\G{\Gamma}
\def\d{\delta}
\def\a{\alpha}
\def\b{\beta}
\def\p{\varphi}
\def\e{\varepsilon}
\def\l{\lambda}
\def\L{\Lambda}
\def\s{\sigma}

\def\o{\omega}

\def\R{{\mathbb R}}
\def\C{{\mathbb C}}
\def\N{{\mathbb N}}
\def\Z{{\mathbb Z}}
\def\Re{\mbox{Re }}
\def\Im{\mbox{Im }}
\def\bs{~\hfill\rule{7pt}{7pt}}

\DeclareMathOperator{\const}{const}
\DeclareMathOperator{\supp}{supp }
\DeclareMathOperator{\dist}{dist}

\newtheorem{Th}{Theorem}
\newtheorem{Pro}{Proposition}
\newtheorem{Def}{Definition}

\newtheorem{Cor}{Corollary}
\begin{document}

\title{Generalized Fourier quasicrystals and almost periodic sets}

\author{Sergii Yu.Favorov}

\address{Sergii Favorov,
\newline\hphantom{iii}  Karazin's Kharkiv National University
\newline\hphantom{iii} Svobody sq., 4, Kharkiv, Ukraine 61022}
\email{sfavorov@gmail.com}

\maketitle {\small
\begin{quote}
\noindent{\bf Abstract.}
Let $\mu$ be a positive measure on the real line with locally finite support $\L$ and integer masses such that its Fourier transform in the sense of distributions is
a purely point measure. An explicit form is found for an entire almost periodic function with a set of zeros $\L$, taking multiplicities into account.
A necessary and sufficient condition for the exponential growth of this function is also found.

    Our constructions are based on the properties of almost periodic sets on the line. In particular, in the article we find a simple representation of such sets.

\medskip


\end{quote}
}

\medskip

   \section{Introduction}\label{S1}
   \bigskip

A crystalline measure on $\R^d$ is a complex measure  with discrete locally finite support, which is temperate distribution and its distributional Fourier transform $\hat\mu$
is also a measure with locally finite support; if, in addition, the measures $|\mu|$ and $|\hat\mu|$ are temperate distributions, then $\mu$ is called the Fourier quasicrystal.

The  Fourier quasicrystal may be considered as a mathematical model for  atomic arrangement
having a discrete diffraction pattern.  There are a lot of papers devoted to study properties of Fourier quasicrystals or, more generally, crystalline measures.
For example, one can mark collections of papers \cite{D}, \cite{Q}, in particular, the basic paper \cite{L1}.

Measures of the form
\begin{equation}\label{a}
\mu=\sum_{\l\in\L} c_\l\d_\l,\qquad c_\l\in\N,
\end{equation}
  are the most important cases of Fourier quasicrystal. Recently A.Olevslskii and A.Ulanovskii \cite{OU1}, \cite{OU2} give a complete description of these measures as zero sets
  of an  exponential polynomial with pure imaginary exponents, where $c_\l$ is just the multiplicity of the zero $\l$ of the polynomial.
 \smallskip

In the present paper we give an analog of the above result for measures \eqref{a} with the distributional Fourier transform
\begin{equation}\label{b}
  \hat\mu=\sum_{\g\in\G}b_\g\d_\g,
\end{equation}
where $\G$ is an arbitrary countable set. In this case the corresponding Poisson's formula
$$
\sum_{\l\in\L} c_\l\hat f(\l)=\sum_{\g\in\G}b_\g f(\g).
$$
also takes place for every function $f$ from Schwarz's class. In order to describe such measures, we use the concept of almost periodic sets, which was introduced
 by M. Krein and B. Levin (\cite{L}, Appendix VI). In modern notations (cf.\cite{M1}, \cite{R}), a locally finite set $\L$ with multiplicities $c_\l$ at points $\l\in\L$
 is almost periodic, if the convolution of measure \eqref{a} with every continuous function with compact support is almost periodic function.
 We will write an almost periodic set like the sequence $A=\{a_n\}_{n\in\Z}$, where each point $a_n=\l$ occurs $c_\l$ times. Therefore, almost periodic sets are in fact multisets.

In Section \ref{S2} we give the original definition of almost periodic sets given by Krein and Levin, which is equivalent to the one given above.
Also we prove some properties of almost periodic sets, in particular, show that such  sets
have the form $\{\a n+\phi(n)\}_{n\in\N}$ with $\a>0$ and an almost periodic mapping $\phi:\,\Z\to\R$.

Note that the zero set of every entire almost periodic function is almost periodic (cf.\cite{T}). On the other hand, it was proved in \cite{FRR}
that every almost periodic set $A\subset\R$ is exactly the zero set of some entire almost periodic function.
In section \ref{S3} we consider measure $\mu$ of the form \eqref{a} under additional conditions that $\hat\mu$ is a measure and $|\hat\mu|$ is a temperate distribution.
If this is the case, the almost periodicity of measure $\mu$, in other words, almost periodicity of corresponding multiset $A$, is equivalent to the property of $\hat\mu$
to be pure point.  Here we find an almost periodic entire function with zero set $A$ in an explicit form depending only on $\g$ and $b_\g$ from equality \eqref{b}.
 Also we find a criterion for $A$ to be the zero set of an almost periodic entire function of the exponential growth.
 Note that according to the Phragmen-Lindelöf principle any entire function bounded on $\R$ cannot grow less than exponentially.

\bigskip

\section{Almost periodic sets}\label{S2}
\bigskip

\begin{Def}[for example, see \cite{C}, \cite{Le}]
 A continuous function on a strip
$$
S=\{z=x+iy:\,-\infty\le a<y<b\le+\infty\}\subset\C
$$
is  almost periodic if for any $\a,\,\b$ such that $[\a,\b]\subset(a,b)$ and
  $\e>0$ the set of $\e$-almost periods
  $$
E_{\a,\b,\e}= \{\tau\in\R:\,\sup_{x\in\R,\a\le y\le\b}|g(x+\tau+iy)-g(x+iy)|<\e\}
  $$
is relatively dense, i.e., $E_{\a,\b,\e}\cap(x,x+L)\neq\emptyset$ for all $x\in\R$ and some $L$ depending on $\e,\a,\b$.
\end{Def}

Just as it was done in \cite{FRR}, we could define an almost periodic set in a strip. But we need it only for the case of sets on the real line, so we can give
this definition here in a simplified form.

\begin{Def}[M. Krein and B. Levin \cite{L}, Appendix VI]\label{D2}
 A discrete locally finite multiset $A=\{a_n\}_{n\in\Z}\subset\R$ is almost periodic if for any $\e>0$ the set of its $\e$-almost periods
\begin{equation}\label{e}
E_\e=\{\tau\in\R:\,\exists\ \text{ a bijection } \s:\Z\to\Z \quad\text{such that}\ \sup_n |a_n+\tau-a_{\s(n)}|<\e\}
\end{equation}
 has nonempty intersection with every interval $(x,x+L_\e)$.
\end{Def}
Set $\mu_A=\sum_n \d_{a_n}$. Clearly, the mass of $\mu_A$ in any point $x\in\R$ is equal to the multiplicity of this point in the sequence $\{a_n\}_{n\in\Z}$.

It was proved in \cite{FRR} that almost periodicity of $A$ is equivalent to almost periodicity of convolution $\mu_A\star\p$ for every $C^\infty$-function $\p(x),\,x\in\R$,
 with compact support. But it is easy to replace $C^\infty$-functions by continuous functions with compact support.

 Indeed, take $C^\infty$-function $\p\ge0$ such that $\p(x)\equiv1$ for $0<x<1$. If $\mu_A\star\p$ is almost periodic, then it is uniformly bounded, hence
 $\mu_A[x,x+1]<K$ for all $x\in\R$ with some constant $K$. For any continuous function $\psi$ with support in $(0,1)$ one can take $\p\in C^\infty$ such that
 $\sup_{x\in\R}|\psi(x)-\p(x)|<\e/K$. We obtain that every $\e$-almost period of $\mu_A\star\p$ is $2\e$-almost period of $\mu_A\star\psi$.

By the way, we gave the proof of the following proposition
\begin{Pro}[\cite{FRR}]\label{P1}
For any almost periodic set there is $k_1\in\N$ such that $\# A\cap[x,x+1]\le k_1$. Also, $\# A\cap[x,x+h)\le k_1(h+1)$.
\end{Pro}
Here and below, $\# H$ means the number of points in the set $H$.
\smallskip

Note that the above proposition implies that $\mu_A(-r,r)=O(r)$ as $r\to\infty$, therefore $\mu_A$ is a temperate distribution.
\begin{Pro}\label{P2}
For any almost periodic set there is $k_2\in\N$ such that for every $h$ and every half-intervals $[x_1,x_1+h),\,[x_2,x_2+h)$ we have
 $|\# A\cap[x_1,x_1+h)-\# A\cap[x_2,x_2+h)|\le k_2$. Also, for every $x\in\R,\, h>0,\,M\in\N$
\begin{equation}\label{m}
|\# A\cap[x,x+h)-(1/M)\# A\cap[x,x+Mh)|\le k_2.
\end{equation}
\end{Pro}
{\bf Proof}. Let $L_1,\,E_1$ be defined in \eqref{e}, and $\tau\in E_1\cap[x_1-x_2,L_1+x_1-x_2)$. Since
$[x_2,x_2+h)+\tau\subset[x_1,x_1+L_1+h)$, we see that to each $a_n\in[x_2,x_2+h)$ assign a point $a_{\s(n)}\in[x_1-1,x_1+L_1+h+1)$. Therefore,
$$
\# A\cap[x_2,x_2+h]\le \# A\cap[x_1,x_1+h)+\# A\cap[x_1-1,x_1)+\# A\cap[x_1+h,x_1+h+L_1+1).
$$
By Proposition \ref{P1}, the last two terms are bounded by $k_1+(L_1+2)k_1$. The proof of the opposite inequality is the same.

To prove the second assertion, we have to add together all the inequalities
$$
\# A\cap[x,x+h)-k_2\le\# A\cap[x+(m-1)h,x+mh)\le k_2+\# A\cap[x,x+h),\quad m=1,2,\dots,M.
$$
\bs

\begin{Pro}\label{P3}
Let $A$ be an almost periodic set. There is a strictly positive density $d$ such that for any $\eta>0$ and any half-interval $I$  with length $l(I)>N_\eta$ we have
$$
            \left|\frac{\# A\cap I}{l(I)}-d\right|<\eta.
$$
\end{Pro}
This result was generalized to all Euclidean spaces in \cite{FK2}.
\medskip

{\bf Proof of Proposition \ref{P3}}. Let $I_1=[x_1,x_1+h_1),I_2=[x_2,x_2+h_2)$ be two half-intervals such that $h_1/h_2=p/q,\,p,q\in\N$. We have
$$
  \frac{\# A\cap I_1}{h_1}-\frac{\# A\cap I_2}{h_2}=\frac{\# A\cap I_1}{h_1}-\frac{\# A\cap qI_1}{qh_1}+\frac{\# A\cap qI_1}{qh_1}-\frac{\# A\cap pI_2}{ph_2}+
  \frac{\# A\cap pI_2}{ph_2}-\frac{\# A\cap I_2}{h_2}.
$$
Applying Proposition \ref{P2}, we get
\begin{equation}\label{i}
 \left|\frac{\# A\cap I_1}{h_1}-\frac{\# A\cap I_2}{h_2}\right|\le \frac{k_2}{h_1}+\frac{k_2}{qh_1}+\frac{k_2}{h_2}\le k_2\left(\frac{2}{h_1}+\frac{1}{h_2}\right).
\end{equation}
Take half-interval $I'=[x_1,x_1+h')$ such that $h_1<h'<h_1+1$ and $h'/h_2$ rational. We have
$$
 \left|\frac{\# A\cap I_1}{h_1}-\frac{\# A\cap I'}{h'}\right|\le\frac{\# A\cap[x_1+h_1,x_1+h')}{h_1}+ \frac{\#A\cap[x_1,x_1+h')}{h_1h'}.
$$
By Proposition \ref{P1}, we obtain
$$
 \left|\frac{\# A\cap I_1}{h_1}-\frac{\# A\cap I'}{h'}\right|\le\frac{k_1}{h_1}+\frac{k_1(h'+1)}{h_1h'}.
$$
Applying \eqref{i} with $I'$ instead of $I_1$, we obtain for all $I_1,\,I_2$
$$
 \left|\frac{\# A\cap I_1}{h_1}-\frac{\# A\cap I_2}{h_2}\right|\le k_2\left(\frac{2}{h_1}+\frac{1}{h_2}\right)+k_1\left(\frac{2}{h_1}+\frac{1}{h_1h'}\right).
$$
Therefore there is a limit
$$
  d=\lim_{l(I)\to\infty}\frac{\# A\cap I}{l(I)}.
$$
Since the set $A$ is relatively dense, this limit is strictly positive.  \bs
\medskip

\begin{Th}\label{T1}
Let $A=\{a_n\}\subset\R$ be an almost periodic set of density $d$ such that $a_n\le a_{n+1}$ for all $n\in\Z$.
Then
\begin{equation}\label{p}
a_n=n/d+\phi(n)\quad\text{with an almost periodic mapping}\quad \phi:\,Z\to\R.
\end{equation}
\end{Th}
{\bf Remark}. The wrong proof of this Theorem was given in \cite{FK1}.
\medskip

{\bf Proof of Theorem \ref{T1}}. After replacing $A$ with $A/d$ we may suppose $d=1$. Also, we may suppose that $a_0<a_1$.
It follows from Proposition \ref{P1} that every interval of length $1$ contains at least one subinterval of length $1/(2k_1)$ that does not intersect $A$.
Take $\e<\min\{1/(6k_1),(a_1-a_0)/3\}$. Divide $\R$ into an infinite number of disjoint semiintervals $I_j=(t_j,t_{j+1}],\,j\in\Z$ such that $t_{j+1}-t_j<2$
and $A\cap (t_j-2\e,t_j+2\e)=\emptyset$ for all $j$.

 Therefore, if $\s$ is any bijection that satisfies \eqref{e}, then $\rho(j)\in\Z$ corresponds to any $j$
such that $\s$ is the bijection of $A\cap I_j$ to $A\cap I_{\rho(j)}$. Let $\s_j$ is the monotone increasing bijection of $A\cap I_j$ on $A\cap I_{\rho(j)}$. Check that
\begin{equation}\label{s0}
|a_n+\tau-a_{\s_j(n)}|<\e\qquad\forall\, a_n\in I_j.
\end{equation}
   Suppose the contrary.  Let $n_0$ be the minimal number such that \eqref{s0} does not satisfy. If $a_{n_0}+\tau+\e\le a_{\s_j(n_0)}$, then
$a_n+\tau+\e\le a_k$ for all $n\le n_0$ and $k\ge\s_j(n_0)$, $a_n\in I_j,\,a_k\in I_{\rho(j)}$. Therefore, $k\neq \s(n)$ for these numbers,
and $\s$ may give  a correspondence
between points from the set $\{n\le n_0:\,a_n\in I_j\}$ and points only from the set $\{k<\s_j(n_0):\,a_k\in I_{\rho(j)}\}$. But
by definition of $\s_j$, we have
$$
\#\{n\le n_0:\,a_n\in I_j\}=\#\{k\le\s_j(n_0):\,a_k\in I_{\rho(j)}\}>\#\{k<\s_j(n_0):\,a_k\in I_{\rho(j)}\}.
$$
We get a contradiction.

If $a_{n_0}+\tau\ge a_{\s_j(n_0)}+\e$, then $a_n+\tau\ge a_k+\e$ for all $n\ge n_0$ and $k\le \s_j(n_0)$, $a_n\in I_j,\,a_k\in I_{\rho(j)}$.
Therefore, $k\neq \s(n)$ for these numbers, and $\s$ may give  a correspondence
between points from the set $\{n\ge n_0:\,a_n\in I_j\}$ and points only from the set $\{k>\s_j(n_0):\,a_k\in I_{\rho(j)}\}$. But
$\#(A\cap I_j)=\#(A\cap I_{\rho(j)})$, hence by definition of $\s_j$, we have
$$
\#\{n\ge n_0:\,a_n\in I_j\}=\#\{k\ge\s_j(n_0):\,a_k\in I_{\rho(j)}\}>\#\{k>\s_j(n_0):\,a_k\in I_{\rho(j)}\}.
$$
We get a contradiction as well.

Since numbers $A\cap I_j$ and $A\cap I_{\rho(j)}$ coincide, we see that the differences between indices of the first elements in these sets coincide for all $j$.
Hence a number $h\in\Z$ corresponds to every $\tau\in E_\e$ such that \eqref{e} satisfies with $\s(n)=n+h$.

It follows from the definition of $\tau$ for all $k\in\N$
$$
   \tau -\e<a_{kh}-a_{(k-1)h}<\tau+\e.
$$
Therefore the length of interval that contains points with numbers from $(k-1)h+1$ till $kh$ is between $\tau-\e$ and $\tau+\e$, and
the length of interval that contains points with numbers from $1$ till $Nh$ is between $N(\tau-\e)$ and $N(\tau+\e)$. Since $d=1$, we get
the inequality $\tau-\e\le h\le\tau+\e$. Set $\phi(n):=a_n-n$. We obtain for all $n\in\Z$
$$
  \phi(n+h)-\phi(n)=a_{n+h}-a_n-h=a_{\s(n)}-(a_n+\tau)+(\tau-h).
$$
Using \eqref{e}, we obtain $|\phi(n+h)-\phi(n)|<2\e$. Therefore, $h$ is $2\e$-almost period of the function $\phi$. The set of $\e$-almost periods $\tau$ of $A$
 is relatively dense, therefore the set of such integers $h$ is relatively dense as well.  \bs

\begin{Cor}\label{C1}
For any almost periodic set $A=\{a_n\}$ such that $0\not\in A$ there is a finite limit
$$
\a_0=\lim_{N\to\infty}\sum_{|a_n|<N}1/a_n.
$$
Moreover, for any $z\in\C\setminus A$ the sum
$$
\a_z=\frac{1}{z-a_0}+\sum_{n\in\N\setminus\{0\}}\left[\frac{1}{z-a_n}+\frac{1}{z-a_{-n}}\right].
$$
converges absolutely.
\end{Cor}
{\bf Proof}. Let $A=\{n/d+\phi(n)\}_{n\in\Z}$. Since the numbers $\phi(n)$ are uniformly bounded, we see that  the sums
$$
 \sum_{n\in\Z,|a_n|<N}\frac{1}{a_n} \quad\mbox{and}\quad \sum_{n\in\Z,|n|<N}\frac{1}{n/d+\phi(n)}
$$
 differ for a uniformly bounded with respect to $N$ number of terms, and each of these terms tends to $0$ as $N\to\infty$. Then
 $$
   \sum_{n\in\Z,0<|n|<N}\frac{1}{n/d+\phi(n)}=\sum_{n\in\N,0<n<N}\frac{\phi(n)+\phi(-n)}{\phi(n)\phi(-n)+n\phi(-n)/d-n\phi(n)/d-(n/d)^2}.
 $$
 The first assertion follows from Cauchy criterium. The second one follows from the absolutely convergence of the series
 $$
 \sum_{n\in\N\setminus\{0\}}\left[\frac{1}{z-a_n}+\frac{1}{z-a_{-n}}\right]=\sum_{n\in\N\setminus\{0\}}\left[\frac{2z-\phi(-n)-\phi(n)}{(n/d+\phi(n)-z)(-n/d+\phi(-n)-z)}\right].
 \phantom{XXX} \bs
 $$
 \medskip

In \cite{L}, Appendix VI, M.Krein and B.Levin considered zero sets $Z_f$ of entire almost periodic functions $f$ of exponential growth. They proved that if
$Z_f\subset\R$, then its zeros $a_n$ form an almost periodic set, which satisfy \eqref{p} and
\begin{equation}\label{p1}
 \sup_{\tau\in\Z}\sum_{n\in\Z\setminus\{0\}}n^{-1}[\phi(n+\tau)-\phi(n)]<\infty.
\end{equation}
On the other hand, they proved that any almost periodic set $A\subset\R$ satisfying conditions \eqref{p} and \eqref{p1} is the set of zeros
of an entire almost periodic function of exponential growth.

It follows from Theorem \ref{T1} that condition \eqref{p} can be omitted in the last result.\medskip

Theorem \ref{T1} was generalized by W.Lawton \cite{La} to almost periodic sets in $\R^d,\,d>1,$, whose spectrum is contained
 in a finitely generated  additive group.

 \section{Almost periodic zeros of entire functions}\label{S3}
\bigskip

 In this section we suppose that $\mu$ is a measure of form \eqref{a}, which is a temperate distribution and its Fourier transform $\hat\mu$ is a measure of form \eqref{b} 
such that $|\hat\mu|$ is also a temperate distribution.
It follows from  \cite {F4}, Lemma 1, that the multiset $A=\{a_n\}_{n\in\Z}$, where each point $a_n=\l$ occurs $c_\l$ times, is an almost periodic set.
In what follows we will suppose that $0\not\in A$. Set
\begin{equation}\label{f}
  f(z)=(1-z/a_0)\prod_{n\in\N} (1-z/a_n)(1-z/a_{-n}).
\end{equation}
It follows from Corollary \ref{C1} that $A$ satisfies Lindelof's condition, hence, $f$ is an entire function of exponential type.

Remark that for any measure $\mu$ of form \eqref{a} with an almost periodic $A$ condition "$\hat\mu$ is a measure" implies "$\hat\mu$ is a pure point measure"
(cf.\cite{M}, Theorem 5.5).
 \smallskip

 Set
 $$
 \R_+:=\{x\in\R:x>0\},\,\R_-:=-\R_+,\,\C_+:=\{z\in\C:\Im z>0\},\,\C_-:=-\C_+,
 $$
 and
$$
a_z(t)=\begin{cases}-2\pi ie^{2\pi itz} &\text{if }t>0,\\0&\text{if }t\le0,\end{cases}\quad z\in\C_+,\quad\qquad
a_z(t)=\begin{cases}2\pi ie^{2\pi itz} &\text{if }t<0,\\0&\text{if }t\ge0,\end{cases}\quad z\in\C_-.
$$
It is not hard to check that in the sense of distributions $\hat a_z(\l)=1/(z-\l)$ for $z\in\C_+\cup\C_-$.

\begin{Pro}\label{P4}
Let $A,\,\mu,\,\hat\mu$ be as above. Then for all $z=x+iy\in\C_+$
\begin{equation}\label{n1}
\frac{f'(z)}{f(z)}=\frac{1}{z-a_0}+\sum_{n\in\N}\left[\frac{1}{z-a_n}+\frac{1}{z-a_{-n}}\right]=-2\pi i\sum_{\g\in\G\cap\R_+}b_\g e^{2\pi i\g z},
\end{equation}
and for all $z=x+iy\in\C_-$
\begin{equation}\label{n2}
\frac{f'(z)}{f(z)}=\frac{1}{z-a_0}+\sum_{n\in\N}\left[\frac{1}{z-a_n}+\frac{1}{z-a_{-n}}\right]=2\pi i\sum_{\g\in\G\cap\R_-}b_\g e^{2\pi i\g z}.
\end{equation}
The function $f'(z)/f(z)$ is almost periodic on each line $y=y_0\neq0$.
\end{Pro}

{\bf Proof}. Let $\p(t)$ be any  even nonnegative $C^\infty$-function such that $\supp\p\subset(-1,1)$ and $\int\p(t)dt=1$. Set $\p_\e(t)=\e^{-1}\p(t/\e)$
for $\e>0$. Fix $z=x+iy\in\C_+$. The functions $a_z(t)\star\p_\e(t)$ and $\hat\a_z(\l)\hat\p_\e(\l)$ belong to Schwartz space of $C^\infty$-functions. Therefore,
\begin{equation}\label{mu}
   (\hat\mu,a_z(t)\star\p_\e(t))=(\mu,\hat\a_z(\l)\hat\p_\e(\l)).
\end{equation}
Then for any  $T_0<\infty$
$$
   (\hat\mu,(a_z\star\p_\e)(t))-(\hat\mu,\a_z(t))
$$
\begin{equation}\label{eps}
 = -2\pi i\sum_{\g\ge T_0}b_\g e^{2\pi i\g z}\int_{-\e}^\e (e^{-2\pi isz}-1)\p_\e(s)ds-2\pi i\sum_{0<\g<T_0}b_\g e^{2\pi i\g z}\int_{-\e}^\e (e^{-2\pi isz}-1)\p_\e(s)ds.
\end{equation}
The first sum is majorized by
\begin{equation}\label{s}
   2\pi(e^{2\pi\e y}+1)\sum_{\g\ge T_0}|b_\g| e^{-2\pi\g y}.
\end{equation}
Set $M(t)=\sum_{\g\in\G:\,0<\g\le t}|b_\g|$. Then
\begin{equation}\label{s1}
   \sum_{\g\ge r}|b_\g|e^{-2\pi\g y}=\int_r^\infty e^{-2\pi ty}dM(t)\le\lim_{T\to\infty}M(T)e^{-2\pi Ty}+2\pi y\int_r^\infty e^{-2\pi ty}M(t)dt.
\end{equation}
 It is easy to check (cf.\cite{F1}) that if $|\hat\mu|$ is a temperate distribution, then
 $|\hat\mu|(-r,r)=O(r^\kappa)$ as $r\to\infty$ with some $\kappa<\infty$. Hence, \eqref{s} is less than any $\eta>0$ for $T_0$ large enough.

 The last sum in \eqref{eps} is less than
  \begin{equation}\label{s2}
 2\pi(e^{2\pi\e y}-1)\sum_{0<\g<T_0}|b_\g|.
 \end{equation}
 Since $\hat\mu$ is a measure, we see that $\sum_{0<\g<T_0}|b_\g|<\infty$, therefore \eqref{s2} is less than $\eta$ for small $\e$. Hence we obtain from \eqref{mu}
 $$
 \lim_{\e\to 0}(\mu,\hat\a_z(\l)\hat\p_\e(\l))=(\hat\mu,\a_z(t))=-2\pi i\sum_{\g\in\G\cap\R_+}b_\g e^{2\pi i\g z}.
$$
 On the other hand, we have
\begin{equation}\label{h}
 (\mu,\hat\a_z(\l)\hat\p_\e(\l))=\frac{\hat\p(\e a_0)}{z-a_0}+\sum_{n\in\N}\left[\frac{\hat\p(\e a_n)}{z-a_n}+\frac{\hat\p(\e a_{-n})}{z-a_{-n}}\right].
\end{equation}
The function $\hat\p(t)$ tends to $1$ as $t\to0$ and $|\hat\p(t)|\le1$.  We have
$$
 \left[\frac{\hat\p(\e a_n)}{z-a_n}+\frac{\hat\p(\e a_{-n})}{z-a_{-n}}\right]=\hat\p(\e a_{-n})\left[\frac{1}{z-a_n}+\frac{1}{z-a_{-n}}\right]+
 \frac{1}{z-a_n}[\hat\p(\e a_n)-\hat\p(\e a_{-n})].
$$
 Since $\hat\p$ is even, we get with bounded $\theta(n)$ and $\phi(n)$
$$
\hat\p(\e a_n)-\hat\p(\e a_{-n})=\hat\p(\e n+\e\phi(n))-\hat\p(-\e n+\e\phi(-n))=\hat\p'(\e n+\e\theta(n))\e|\phi(n)-\phi(-n)|.
$$
Since $\hat\p(t)$ belongs to Schwartz space,  we see that $\hat\p'(t)=O(1/|t|)$ as $t\to\infty$. Hence for $\e\ge 1/|n+\theta(n)|$
$$
|\e[\hat\p'(\e n+\e\theta(n))]|\le C|n|^{-1}
$$
with a constant $C<\infty$. The same estimate (with another constant $C$) is valid for $\e<1/|n+\theta(n)|$,
that for all $n\in\N$ and $\e>0$
$$
|\hat\p(\e a_n)-\hat\p(\e a_{-n})|\le (C/n)2\sup_n|\phi(n)|.
$$
 Hence the right-hand side of \eqref{h} for all $\e>0$ is majorized by the  sum
\begin{equation}\label{sum}
\frac{1}{|z-a_0|}+\sum_{n\in\N}\left|\frac{1}{z-a_n}+\frac{1}{z-a_{-n}}\right|+\sum_{n\in\N}\frac{C'}{n|z-a_n|}.
\end{equation}
By Theorem \ref{T1}, we have  $1/(z-a_n)=O(1/n)$. Taking into account also Corollary \ref{C1}, we get the convergence of both sums in \eqref{sum}.
Therefore we can go to the limit  in \eqref{h} as $\e\to0$ and obtain \eqref{n1}.

 By \eqref{s1}, $\sum_{\g\ge 1}|b_\g|e^{-2\pi\g y_0}<\infty$ for $y_0>0$,
and $\sum_{0<\g<1}|b_\g|<\infty$. Therefore the series in right-hand part of \eqref{n1} absolutely converges, and $f'(z)/f(z)$ is almost periodic on the line $y=y_0$.

Furthermore,  the measure $\mu$ is real-valued, hence, $b_{-\g}=\bar b_\g$.
Therefore, in the case $y_0<0$ we can apply \eqref{n1} to the function $\overline{f(\bar z)}$ and obtain \eqref{n2}. \bs

\begin{Th}\label{T2}
For measures $\mu$, $\hat\mu$, and the almost periodic set $A=\{a_n\}_{n\in\Z}$ as above there is the  almost periodic entire function of the form
\begin{equation}\label{F}
    F(z)=e^{g(z)}f(z),
\end{equation}
where $f(z)$ is defined in \eqref{f} and
\begin{equation}\label{g}
g(z):=\sum_{\g\in\G,0<\g<1}b_\g\frac{e^{2\pi i\g z}-1}{\g}.
\end{equation}
\end{Th}

{\bf Proof}. The sum  $\sum_{\g\in\G\cap\R_+}b_\g e^{2\pi i\g z}$ is majorized  by
$$
 \sum_{\g\in\G,o<\g<1}|b_\g|+\sum_{\g\ge1}|b_\g|e^{-2\pi\g y}.
$$
Since $\hat\mu$ is a measure, the first sum is finite, the second one is also finite due to \eqref{s1}. In particular,  the function
$$
  g(z)=\sum_{k=1}^\infty\frac{(2\pi iz)^k}{k!}\sum_{\g\in\G,0<\g<1}\g^{k-1}b_\g
$$
 is well-defined and holomorphic for all $z\in\C$. Using \eqref{n1} for $z=x+iy\in\C_+$ and changing the order of summing and integrating,  we get
\begin{equation}\label{pr}
  \log f(z)-\log f(iy)=\int_0^x\frac{f'(t+iy)}{f(t+iy)}dt=-\left[\sum_{\g\in\G\cap\R_+}b_\g e^{-2\pi\g y}\frac{e^{2\pi i\g x}-1}{\g}\right].
\end{equation}
Also, both sums
$$
        \sum_{\g\in\G,\g\ge1}\g^{-1}e^{-2\pi\g y}\,b_\g(1-e^{2\pi i\g x})\quad\mbox{and}\quad\sum_{\g\in\G,0<\g<1}\frac{e^{-2\pi\g y}-1}{\g}b_\g(1-e^{2\pi i\g x})
$$
are uniformly bounded in every strip $\{z:\,0<\a\le\Im z\le\b<\infty\}$. Therefore in this strip
\begin{equation}\label{c}
   \log f(z)=-\sum_{\g\in\G,0<\g<1}b_\g \g^{-1}[e^{2\pi i\g z}-1]+O(1).
\end{equation}
In particular, the function $\log|f(z)|+\Re g(z)$ is bounded in this strip, and the function $F(z)$ is uniformly bounded in every closed
horizontal substrip of a finite width in $\C_+$.

Arguing by the same way, we get for $z=x+iy\in\C_-$
$$
   \log f(z)=\sum_{\g\in\G,-1<\g<0}b_\g \g^{-1}[e^{2\pi i\g z}-1]+O(1).
$$
Since $b_{-\g}=\bar b_\g$, we see that the real part of the right-hand side of the equality equals
$$
-\Re\sum_{\g\in\G,0<\g<1}\bar b_\g \g^{-1}[e^{-2\pi i\g z}-1]+O(1)=-\Re\sum_{\g\in\G,0<\g<1}b_\g \g^{-1}[e^{2\pi i\g\bar z}-1]+O(1).
$$
Taking into account that
$$
  \left|\sum_{\g\in\G,0<\g<1}b_\g \g^{-1}[e^{2\pi i\g\bar z}-1]-\sum_{\g\in\G,0<\g<1}b_\g \g^{-1}[e^{2\pi i\g z}-1]\right|\le
  \sum_{\g\in\G,0<\g<1}|b_\g|\left|\frac{e^{2\pi\g y}-e^{-2\pi\g y}}{\g y}\right||y|,
$$
we obtain that the function $\log|f(z)|+\Re g(z)$ is bounded in every strip $\{z:\,-\infty<\a\le\Im z\le\b<0\}$,
and the function $F(z)$ is uniformly bounded in every closed horizontal substrip of a finite width in $\C_-$ as well.

Furthermore, in any strip $\{z:\,|\Im z|<M\}$
\begin{equation}\label{g1}
  |g(z)-g(x)|\le  \sum_{\g\in\G,0<\g<1}|b_\g|\left|\frac{e^{2\pi i\g x}(e^{2\pi\g y}-1)}{\g y}\right||y|<C(M),
\end{equation}
and
\begin{equation}\label{g2}
   |g(x)|\le\left[\sum_{\g\in\G,0<\g<\e}+\sum_{\g\in\G,\e\le\g<1}\right]|b_\g|\left|\frac{1- e^{2\pi i\g x}}{\g}\right|
\le 2\pi|x|\sum_{o<\g<\e}|b_\g|+2\e^{-1}\sum_{\e\le\g<T_0}|b_\g|.
\end{equation}
Since $\sum_{o<\g<\e}|b_\g|$ is arbitrary small, we see that $g(x)=o(|x|)$ as $x\to\infty$. Therefore the function $g(z)$ has the same growth
in any horizontal strip of the finite width. Also, $\log|f(z)|\le O(|z|)$ in this strip.
Applying Fragment-Lindelof Principle, we obtain that  $F(z)$ is bounded in every such strip.

It follows from \eqref{c} that the function $\log F(z)=\log f(z)+g(z)$ is bounded on any line $y=\const$. Then we have
$$
(\log F(x+iy))' =f'(x+iy)/f(x+iy)+2\pi i\sum_{\g\in\G,0<\g<1}b_\g e^{2\pi i\g x-2\pi y}.
$$
Since the function  $f'(x+iy)/f(x+iy)$ is almost periodic in $x$, we see that Bohr's theorem (cf.\cite{Le}, Theorem 1.2.1) implies almost periodicity of the function $\log F(x+iy)$
and then $F(x+iy)$.
By (\cite{Le}, Theorem 1.2.3), the function $F(z)$ is almost periodic in every strip, where it is bounded,
hence it is almost periodic in $\C$. \bs

\begin{Th}\label{T3}
In conditions of the previous theorem, $A$ to be the zero set of an almost periodic entire function  of exponential type if and only if the function $g(z)$ from \eqref{g}
is uniformly bounded on $\R$. If this is the case, the entire function \eqref{f} is almost periodic.
\end{Th}
{\bf Proof}. Let $g(z)$ be uniformly bounded for $z\in\R$. By \eqref{g1}, the function $g(x+iy)$ is also bounded in $x$ for any fixed $y>0$.
It follows from \eqref{c} that the function $\log f(x+iy)$ is also bounded in $x\in\R$. But $(\log f(x+iy))'$ is almost periodic,
hence, by Bohr's Theorem, the functions $\log f(x+iy)$ and $f(x+iy)$ are almost periodic too. The function $f$ is the exponential type,
and by Fragment-Lindelof Principle, it is bounded on every horizontal strip of a finite width. Hence, it is almost periodic on such strips, consequently, $f$
is an almost periodic entire function with the given zero set $A$.

 Now suppose that  $G(z)$ is an entire almost periodic function of exponential type with zero set $A$. Clearly, $G(z)=K_1e^{K_2z}f(z)$ with $K_1,\,K_2\in\C$.
  Since zero set of $G(z)$ coincides with $A$,
 we get, using Lemma 1 from \cite{L}, Ch.6, that for every $\e>0$ and $M<\infty$ there is $m(\e)>0$ such that
 $$
 |G(z)|\ge m(\e)\quad\mbox{for}\quad \{z:\,|\Im z|\le M,\,z\not\in A(\e)\},\quad\mbox{where}\quad A(\e):=\{z:\,\dist(z,A)<\e\}.
 $$
  By Proposition \ref{P1}, for $\e$ small enough each connected component of $A(\e)$ contains no segment of length $1$, hence its diameter is less than $1$.
 Let $F$ be the almost periodic function that defined in \eqref{F}. The holomorphic function $F(z)/G(z)$ is uniformly bounded on the set $\{z:\,|\Im z|<M,\,z\not\in A(\e)\}$, therefore
 it is uniformly bounded on the whole strip $|\Im z|<M$.
Then $|G(x+iy)|\ge m(\e)$ for any $y=y_0>\e$, therefore $F(x+iy)/G(x+iy)$ is almost periodic in $x$ for $y=y_0$. Consequently, it is almost periodic for every $|y|<M$,
 in particular, on the real line. Moreover, $F(z)/G(z)$ has no zeros, hence the same Lemma 1 from \cite{L} implies that $|F(x)/G(x)|\ge c>0$ for all $x$.
 Now by Theorem 2.7.1 from \cite{Le},
 $$
 F(x)/G(x)=e^{h(x)+i\o x},\qquad \o\in\R,
 $$
 with almost periodic $h(x)$. Therefore the function $g(x)-K_2x-i\o x=h(x)$ is bounded on $\R$. By \eqref{g2},  $g(x)=o(|x|)$ as $x\to\infty$, hence,
 $K_2=-i\o$ and $g(x)$ is bounded on $\R$. \bs

\begin{Cor}\label{C2}
If $\sum_{\g\in\G,0<\g<1}\g^{-1}|b_\g|<\infty$, then $A$ is the zero set of some almost periodic entire function  of exponential type.
\end{Cor}

\end{document}